\documentclass[12pt]{article}
\usepackage{mathrsfs}
\usepackage{amssymb,a4}
\usepackage{amsmath,amsfonts,amssymb,amsthm}
\newcommand{\h}{\mathfrak{h}}
\newcommand{\al}{\alpha}

\newcommand{\1}{{{\bf 1}}}

\newcommand{\End}{{\rm End}\,}

\def\wt{{\rm wt}}

\def\de{\delta}
\def\be{\beta}

\def\C{{\mathbb C}}
\def\Q{{\mathbb Q}}
\def\c{{\tilde{c}}}

\def\Z{{\mathbb Z}}

\def\N{{\mathbb N}}
\def\1{{\bf 1}}
\def\l{\lambda}
\def \wt{{\rm wt}}

\def \End{{\rm End}}

\def \<{\langle}
\def \>{\rangle}
\def \w{\omega}
\def \pf{\noindent {\bf Proof: \,}}
\def\theequation{5.\arabic{equation}}

\renewcommand{\theequation}{\thesection.\arabic{equation}}

\newtheorem{theorem}{Theorem}[section]
\newtheorem{prop}[theorem]{Proposition}
\newtheorem{lem}[theorem]{Lemma}
\newtheorem{coro}[theorem]{Corollary}
\newtheorem{remark}[theorem]{Remark}
\theoremstyle{definition}
\newtheorem{definition}[theorem]{Definition}

\begin{document}
\begin{center}
{\Large {\bf A characterization of the vertex operator
algebra $V_{L_{2}}^{A_{4}}$}} \\
\vspace{0.5cm}
 Chongying Dong \footnote{Supported bf NSF grants}\\
 Department of Mathematics,  University of California, Santa Cruz, CA 95064, USA\\
 Cuipo Jiang\footnote{Supported  by China NSF grants
10931006,11371245, the RFDP of China(20100073110052), and the
Innovation Program of Shanghai Municipal
Education Commission (11ZZ18)}\\
Department of Mathematics, Shanghai Jiaotong University, Shanghai 200240 China

%\textcolor[rgb]{0.50,0.50,1.00}{dong@math.ucsc.edu; cpjiang@sjtu.edu.cn}}
\end{center}
\hspace{1cm}
\begin{abstract} The rational vertex operator algebra $V_{L_{2}}^{A_{4}}$ is characterized in terms of  weights of primary vectors. This reduces the classification of rational vertex operator algebras with $c=1$ to the characterizations of $V_{L_{2}}^{S_{4}}$ and $V_{L_{2}}^{A_{5}}.$

2000MSC:17B69
\end{abstract}

\section{Introduction}

Characterizations of vertex operator algebras $V_{L_2}^G$ for root lattice $L_2$ of $sl(2,\C)$ and
finite groups $G=A_4, S_4, A_5$  are the remaining part of classification of rational vertex operator algebras with $c=1$ after the work of \cite{DM2}, \cite{ZD}, \cite{DJ1}-\cite{DJ3}. Using the structure and representation theory of $V_{L_2}^{A_4}$ obtained in \cite{DJ4} and \cite{DJJJY}, we give a characterization of rational vertex operator algebra $V_{L_2}^{A_4}$ in this paper.

The main assumption for the  characterization of vertex operator algebra $V_{\Z\gamma}^+$ with $(\gamma,\gamma)\geq 6$ being a positive even integer in \cite{DJ2}-\cite{DJ3} is that the dimension of the weight 4 subspace is at least three dimensional. Knowing  the explicit structure of $V_{L_2}^{A_4}$ we have a different assumption in
characterizing $V_{L_2}^{A_4}.$ That is, there is a primary vector of weight $9$ and the weight of any primary vector which is not a multiple of $\1$  is greater than or equal to $9.$ Due to a recent result in \cite{DLN} on the modularity of the $q$-characters of the irreducible modules for rational and $C_2$-cofinite vertex operator algebras, we can use the classification of $q$-characters of rational vertex operator algebras with $c=1$ from \cite{K} to conclude that the $q$-character of such a vertex operator algebra and that of $V_{L_2}^{A_4}$ are the same.

Two basic facts are used in the characterization. The first one is that both $V_{L_2}^{A_4}$
and an abstract vertex operator algebra $V$ satisfying the required conditions have the same decomposition as modules for the Virasoro algebra. This allows us to use the fusion rules for the irreducible modules
for the Virasoro vertex operator algebra $L(1,0)$ obtained in \cite{DJ1} and \cite{M} to understand the structure of both vertex operator algebras in terms of generators. The other one is that $V_{L_2}^{A_4}$
is generated by a weight $9$ primary vector $x^1$ and has a spanning set in terms Virasoro algebra, the component operators of $x^1$ and the component operators of $y^1$ which is a primary vector of weight $16.$
From the $q$-characters we know that $V$ also has primary vectors $x^2,y^2$ of weights $9,$ $16,$ respectively. The main task is to show how the vertex operator algebra $V$ has a similar spanning set with $x^1,y^1$ being replaced by $x^2,y^2.$  The fusion rules for the vertex operator algebra $L(1,0)$ play a crucial role here.

We certainly expect that the ideas and methods presented in this paper can also be used to characterize vertex
operator algebras $V_{L_2}^G$ for $G=S_4, A_5$ although $V_{L_2}^G$ have not been  understood well. It seems that knowing the generators of $V_{L_2}^G$ and a spanning set is good enough for the purpose of characterization. Of course, the rationality is also needed.

The paper is organized as follows.  We review the modular invariance results from \cite{Z} and \cite{DLN} in Section 2. These results will be used to conclude that $V_{L_2}^{A_4}$ and an abstract vertex operator algebra $V$ satisfying certain conditions  have the same graded dimensions. We also review the fusion rules for
the vertex operator algebra $L(1,0)$ from \cite{M} and \cite{DJ1} in this section.  In Sections 3 we discuss the structure of $V_{L_2}^{A_4}$ including the generators
and spanning set following \cite{DJ4}. Section 4 is devoted to the characterization of $V_{L_2}^{A_4}.$ That is, if a rational, $C_2$-cofinite and self-dual vertex operator algebra $V$ of central charge 1 satisfies (a)
$V$ is a completely reducible module for the Virasoro algebra, (b) $V$ has a primary vector of weight $9$ and
the weight of  any primary vector whose weight is greater than $0$ is greater than or equal to $9,$ then $V$ is isomorphic to  $V_{L_2}^{A_4}.$
The main idea is to use generators and relations to construct a vertex operator algebra isomorphism from
 $V_{L_2}^{A_4}$ to $V.$

\section{Preliminaries}
\def\theequation{2.\arabic{equation}}
\setcounter{equation}{0}

Let $V=(V,Y,{\bf 1},\omega)$ be a vertex operator algebra \cite{B},
\cite{FLM}. We review various notions of $V$-modules (cf.
\cite{FLM}, \cite{Z}, \cite{DLM1}) and the definition of  rational
vertex operator algebras.  We also discuss some consequences
following \cite{DLM2}, \cite{K}, \cite{DLN}, \cite{DJ2}.

\begin{definition} A weak $V$ module is a vector space $M$ equipped
with a linear map
$$
\begin{array}{ll}
Y_M: & V \rightarrow {\rm End}(M)[[z,z^{-1}]]\\
 & v \mapsto Y_M(v,z)=\sum_{n \in \Z}v_n z^{-n-1},\ \ v_n \in {\rm End}(M)
\end{array}
$$
satisfying the following:

1) $v_nw=0$ for $n>>0$ where $v \in V$ and $w \in M$

2) $Y_M( {\textbf 1},z)=Id_M$

3) The Jacobi identity holds:
\begin{eqnarray}
& &z_0^{-1}\de \left(\frac{z_1 - z_2}{
z_0}\right)Y_M(u,z_1)Y_M(v,z_2)-
z_0^{-1} \de \left(\frac{z_2- z_1}{ -z_0}\right)Y_M(v,z_2)Y_M(u,z_1) \nonumber \\
& &\ \ \ \ \ \ \ \ \ \ =z_2^{-1} \de \left(\frac{z_1- z_0}{
z_2}\right)Y_M(Y(u,z_0)v,z_2).
\end{eqnarray}
\end{definition}

%admissible

\begin{definition}
An admissible $V$ module is a weak $V$ module  which carries a
$\Z_+$-grading $M=\bigoplus_{n \in \Z_+} M(n)$, such that if $v \in
V_r$ then $v_m M(n) \subseteq M(n+r-m-1).$
\end{definition}

\begin{definition}
An ordinary $V$ module is a weak $V$ module which carries a
$\C$-grading $M=\bigoplus_{\l \in \C} M_{\l}$, such that:

1) $dim(M_{\l})< \infty,$

2) $M_{\l+n}=0$ for fixed $\l$ and $n<<0,$

3) $L(0)w=\l w=\wt(w) w$ for $w \in M_{\l}$ where $L(0)$ is the
component operator of $Y_M(\omega,z)=\sum_{n\in\Z}L(n)z^{-n-2}.$
\end{definition}

\begin{remark} \ It is easy to see that an ordinary $V$-module is an admissible one. If $W$  is an
ordinary $V$-module, we simply call $W$ a $V$-module.
\end{remark}

We call a vertex operator algebra rational if the admissible module
category is semisimple. We have the following result from
\cite{DLM2} (also see \cite{Z}).

\begin{theorem}\label{tt2.1}
If $V$ is a  rational vertex operator algebra, then $V$ has finitely
many irreducible admissible modules up to isomorphism and every
irreducible admissible $V$-module is ordinary.
\end{theorem}

Suppose that $V$ is a rational vertex operator algebra and let
$M^1,...,M^k$ be the irreducible  modules such that
$$M^i=\oplus_{n\geq 0}M^i_{\l_i+n}$$
where $\l_i\in\Q$ \cite{DLM2},  $M^i_{\l_i}\ne 0$ and each
$M^i_{\l_i+n}$ is finite dimensional. Let $\l_{min}$ be the minimum
of $\l_i$'s. The effective central charge $\c$ is defined as
$c-24\l_{min}.$ For each $M^i$ we define the $q$-character of $M^i$
by
$$Z_{i}(q)=q^{-c/24}\sum_{n\geq 0}(\dim M^i_{\l_i+n})q^{n+\l_i}.$$

A vertex operator algebra $V$ is called $C_2$-cofinite if $\dim V/C_2(V)$ is finite dimensional where
$C_2(V)$ is a subspace of $V$ spanned by $u_{-2}v$ for $u,v\in V.$
If $V$ is $C_2$-cofinite, then $Z_{i}(q)$ converges to
a holomorphic function on $0<|q|<1$ \cite{Z}. Let $q=e^{2\pi i\tau}$ and we sometimes also write $Z_i(q)$ by $Z_i(\tau)$ to indicate that $Z_i(q)$ is a holomorphic function on the upper half plane.

For a $V$-module $W$, let
$W^{\prime}$ denote the graded dual of $W$. Then $W'$ is also a
$V$-module \cite{FHL}. A vertex operator algebra $V$ is called self dual if $V'$ \cite{FHL} is isomorphic to itself.
The following result comes from \cite{DLN}
\begin{theorem}\label{dln} Let $V$ be a rational, $C_2$-cofinite, self dual simple vertex operator algebra.

(1) Each $Z_i(\tau)$ is a modular function on a congruence subgroup of $SL_2(\Z)$ of level $n$ which is the smallest positive integer such that $n(\lambda_i-c/24)$ is an integer for all $i.$

(2) $\sum_{i}|Z_i(\tau)|^2$ is $SL_2(\Z)$-invariant.
\end{theorem}

We now recall the construction of vertex operator algebras $M(1)^+,$ $V_L^+$ and
related results from \cite{A1}, \cite{A2}, \cite{AD}, \cite{ADL},
\cite{DN1}, \cite{DN2}, \cite{DN3}, \cite{DJL},
 \cite{FLM}.

Let $L=\Z \alpha$ be a positive definite lattice with $(\alpha,\alpha)=2k$ for some positive integer $k.$  Set
$\h=\C\otimes_{\Z} L$ and extend $(\cdot\,,\cdot)$ to $\h$ by
$\C$-linearity. Let
$\hat{\h}=\C[t,t^{-1}]\otimes\h\oplus\C K$ be the corresponding affine Lie algebra
so that
\begin{align*}
[\alpha(m),\,\alpha(n)]=2km\delta_{m+n,0}K\hbox{ and }[K,\hat{\h}]=0
\end{align*}
for any $m,\,n\in\Z$, where $\alpha(m)=\alpha\otimes t^m.$  Note that
$\hat{\h}^{\geq 0}=\C[t]\otimes\h\oplus\C K$ is an abelian
subalgebra. Let $\C e^\lambda$ (for any $\lambda\in\h$) be one-dimensional
$\hat{\h}^{\geq 0}$-module such that $\alpha(m)\cdot
e^\lambda=(\lambda,\alpha)\delta_{m,0}e^\lambda$ and $K\cdot
e^\lambda=e^\lambda$ for $m\geq0$. Consider the induced module
\begin{align*}
M(1,{\lambda})=U(\hat{\h})\otimes_{U(\hat{\h}^{\geq 0})}\C
e^\lambda\cong S(t^{-1}\C[t^{-1}])\ ({\rm linearly}).
\end{align*}
 Set
$$M(1)=M(1,0).$$ Then there exists a linear map
$Y:M(1)\to\End M(1)[[z,z^{-1}]]$ such that $(M(1),\,Y,\,\1,\,\w)$
is a simple vertex operator algebra and
$M(1,\lambda)$ is an irreducible $M(1)$-module for any
$\lambda\in\h$ (see \cite{FLM}). The vacuum vector and the Virasoro
element are given by $\1=e^0$ and $\w=\frac{1}{4k}\alpha(-1)^2\1,$
respectively.

We use  $\C[L]$ to denote  the group algebra of $L$ with a basis $e^{\beta}$ for
$\beta\in L.$ Then
$$V_L=M(1)\otimes \C[L]$$
is the lattice vertex operator algebra associated to $L$ \cite{B}, \cite{FLM}.
Let $L^{\circ}$ be the dual lattice  of $L:$
$$L^{\circ}=\{\,\lambda\in\h\,|\,(\alpha,\lambda)\in\Z\,\}=\frac{1}{2k}L$$
and $L^{\circ}=\cup_{i=-k+1}^k(L+\lambda_i)$ be the coset
decomposition with $\lambda_i=\frac{i}{2k}\alpha.$ Set $\C[L+\lambda_i]=\bigoplus_{\beta\in L}\C
e^{\beta+\lambda_i}.$ Then each $\C[L+\lambda_i]$ is an
$L$-submodule in an obvious way. Set
$V_{L+\lambda_i}=M(1)\otimes\C[L+\lambda_i]$. Then $V_L$ is a
rational vertex operator algebra and $V_{L+\lambda_i}$ for
$i=-k+1,\cdots,k$ are the irreducible modules for $V_L$ (see
\cite{B}, \cite{FLM}, \cite{D1}).

Let $\theta:V_{L+\lambda_i}\to
V_{L-\lambda_i}$ be a linear isomorphism for $i\in\{-k+1,\cdots,k\}$ such that
\begin{align*}
\theta(\alpha(-n_{1})\alpha(-n_{2})\cdots \alpha(-n_{s})\otimes
e^{\beta+\lambda_i})=(-1)^{s}\alpha(-n_{1})\alpha(-n_{2})\cdots
\alpha(-n_{s})\otimes e^{-\beta-\lambda_i}
\end{align*}
where $n_j>0$ and $\beta\in L.$   In particular, $\theta$
is an automorphism of $V_{L}$ which induces an automorphism of
$M(1).$ For any $\theta$-stable subspace $U$ of $V_{L^{\circ}}$, let $U^\pm$
be the $\pm1$-eigenspace of $U$ for $\theta$. Then $V_L^+$ is a
simple vertex operator algebra.

The $\theta$-twisted Heisenberg algebra $\h[-1]$ and its
irreducible module $M(1)(\theta)$ from \cite{FLM} are also needed. Define  a
character $\chi_s$ of $L/2L$ such that $\chi_s(\alpha)=(-1)^s$ for $s=0,1$
and let $T_{\chi_s}=\C$ be the corresponding irreducible $L/2L.$ Then
$V_L^{T_{\chi_s}}=M(1)(\theta)\otimes T_{\chi_s}$ is an irreducible
$\theta$-twisted $V_L$-module (see \cite{FLM}, \cite{D2}). We also define
actions of $\theta$ on  $M(1)(\theta)$ and $V_L^{T_{\chi_s}}$ by
\begin{align*}
\theta(\alpha(-n_{1})\alpha(-n_{2})\cdots
\alpha(-n_{p}))=(-1)^{p}\alpha(-n_{1})\alpha(-n_{2})\cdots
\alpha(-n_{p})
\end{align*}
\begin{align*}
\theta(\alpha(-n_{1})\alpha(-n_{2})\cdots \alpha(-n_{p})\otimes
t)=(-1)^{p}\alpha(-n_{1})\alpha(-n_{2})\cdots \alpha(-n_{p})\otimes
t
\end{align*}
for $n_j\in \frac{1}{2}+\Z_{+}$ and $t\in T_{\chi_s}$. We denote the
$\pm 1$-eigenspaces of $M(1)(\theta)$ and $V_L^{T_{\chi_s}}$ under
$\theta$ by $M(1)(\theta)^{\pm}$ and $(V_L^{T_{\chi_s}})^{\pm}$
respectively.

The classification of irreducible modules for arbitrary $M(1)^+$ and $V_L^+$ are obtained in \cite{DN1}-\cite{DN3}
and \cite{AD}. The rationality of $V_L^+$ is established in
\cite{A2} for rank one lattice  and \cite{DJL} in general. One can find the following results from these papers.
\begin{theorem}\label{t32} (1)
Any irreducible module for the vertex operator algebra $M(1)^+$ is
isomorphic to one of the following modules:$$ M(1)^+, M(1)^-, M(1,
\lambda) \cong M(1, -\lambda)\ (0\neq \lambda \in \h),
M(1)(\theta)^+, M(1)(\theta)^- .$$

(2) Any irreducible $V_L^+$-module is isomorphic to one of the following
modules:
$$V_L^{\pm}, V_{\lambda_i+L}( i \not= k),
V_{\lambda_k+L}^{\pm}, (V_L^{T_{\chi_s}})^{\pm}.$$

(3)  $V_{L}^{+}$ is rational.
\end{theorem}

 The following characterization of $V_{L}^+$  is given in \cite{DJ2} and \cite{DJ3}.
\begin{theorem}
Let $V$ be a simple, rational  and $C_{2}$-cofinite self-dual vertex operator
algebra such that $V$ is generated by highest vectors of the
Virasoro algebra $L(1,0)$ with $\tilde{c}=c=1$ and
$$\dim V_{2}=1, \ \dim V_{3}\geq 2$$
or
$$
\dim V_{2}=\dim V_{3}=1, \  \dim V_{4}\geq 3.$$ Then $V$ is
isomorphic to $V_{\Z\al}^{+}$ for some rank one positive definite
even lattice $L=\Z\al$.
\end{theorem}

We will need the following result from \cite{K}.
\begin{theorem}\label{Kir}
 Let $V$ be a rational CFT type vertex operator algebra
with $c=\c=1$ such that each $Z_i(\tau)$ is a modular function on a congruence subgroup and
$\sum_i |Z_{i}(\tau)|^2$ is $SL_2(\Z)$-invariant, then the q-character of $V$ is equal to the character of
one of the following vertex operator algebras $V_{L}, V_{L}^{+}$ and
$V_{{\mathbb
Z}\alpha}^{G}$,
 where $L$ is any positive definite even lattice of
rank 1, ${\mathbb Z}\alpha$ is the root lattice of type $A_{1}$ and
$G$ is a finite subgroup of $SO(3)$ isomorphic to $A_{4}, S_{4}$ or
$A_{5}$.
\end{theorem}

By Theorems \ref{dln} we know that the assumptions in Theorem \ref{Kir} hold. So the $q$-character of a rational vertex operator algebra with $c=1$ is known.

Recall from \cite{FHL} the fusion rules of
vertex operator algebras. Let $V$ be a vertex operator algebra, and
$ W^i$ $ (i=1,2,3$)
 be  ordinary $V$-modules. We denote by $I_{V} \left(\hspace{-3 pt}\begin{array}{c} W^3\\
W^1\,W^2\end{array}\hspace{-3 pt}\right)$  the vector space of all
intertwining operators of type $\left(\hspace{-3 pt}\begin{array}{c}
W^3\\ W^1\,W^2\end{array}\hspace{-3 pt}\right)$.
  It is well known that fusion rules have the
following symmetry \cite{FHL}.

\begin{prop}\label{p4.2}
Let $W^{i}$ $(i=1,2,3)$ be $V$-modules. Then
$$\dim I_{{V}} \left(\hspace{-3 pt}\begin{array}{c} W^3\\
W^1\,W^2\end{array}\hspace{-3 pt}\right)=\dim I_{{V}} \left(\hspace{-3 pt}\begin{array}{c} W^3\\
W^2\,W^1\end{array}\hspace{-3 pt}\right), \ \ \ \dim I_{{V}} \left(\hspace{-3 pt}\begin{array}{c} W^3\\
W^1\,W^2\end{array}\hspace{-3 pt}\right)=\dim I_{{V}} \left(\hspace{-3 pt}\begin{array}{c} (W^2)^{\prime}\\
W^1\,(W^3)^{\prime}\end{array}\hspace{-3 pt}\right).$$
\end{prop}

Here are some results on the fusion rules for the Virasoro vertex operator algebra. Recall that $L(c,h)$ is the irreducible highest weight module for
the Virasoro algebra with central charge $c$ and highest weight $h$
for $c,h\in \C.$ It is well known that $L(c,0)$ is a vertex operator
algebra \cite{FZ}. The following two results can be found  in \cite{M} and
\cite{DJ1}.
\begin{theorem}\label{2t1} (1) We have
$$
\dim I_{L(1,0)} \left(\hspace{-3 pt}\begin{array}{c} L(1,k^{2})\\
L(1, m^{2})\,L(1, n^{2})\end{array}\hspace{-3 pt}\right)=1,\ \
k\in{\mathbb Z}_{+},  \ |n-m|\leq k\leq n+m,$$
$$\dim I_{L(1,0)} \left(\hspace{-3 pt}\begin{array}{c} L(1,k^{2})\\
L(1, m^{2})\,L(1, n^{2})\end{array}\hspace{-3 pt}\right)=0,\ \
k\in{\mathbb Z}_{+},  \ k<|n-m| \ {\rm or} \  k>n+m, $$ where
$n,m\in{\mathbb Z}_{+}$.

(2) For $n\in{\mathbb Z}_{+}$ such that $n\neq p^{2}$, for all
$p\in{\mathbb Z}_{+}$, we have
$$\dim I_{L(1,0)} \left(\hspace{-3 pt}\begin{array}{c} L(1,n)\\
L(1, m^{2})\,L(1, n)\end{array}\hspace{-3 pt}\right)=1,$$
$$\dim I_{L(1,0)} \left(\hspace{-3 pt}\begin{array}{c} L(1,k)\\
L(1, m^{2})\,L(1, n)\end{array}\hspace{-3 pt}\right)=0,$$
 for $k\in{\mathbb Z}_{+}$ such that  $k\neq n$.

\end{theorem}

\section{The vertex operator subalgebra $V_{L_{2}}^{A_{4}}$}
\def\theequation{3.\arabic{equation}}
\setcounter{equation}{0} Let $L_{2}=\Z\al$ be the rank one
positive-definite lattice such that $(\al,\al)=2.$
Then $(V_{L_{2}})_{1}$ is a Lie algebra isomorphic to $sl_{2}(\C)$ and  has an
orthonormal basis:
$$x^{1}=\frac{1}{\sqrt{2}}\al(-1){\bf 1}, \ x^2= \frac{1}{\sqrt{2}}(e^{\al}+e^{-\al}), \
x^3=\frac{i}{\sqrt{2}}(e^{\al}-e^{-\al}).$$
There are three involutions $\tau_{i}\in
Aut(V_{L_{2}})$, $i=1,2,3$ be such that
$$
\tau_{1}(x^1,x^2,x^3)=(x^1,x^2,x^3)\left[\begin{array}{ccc}1& &\\
&-1&\\&&-1\end{array}\right],$$
$$
\tau_{2}(x^1,x^2,x^3)=(x^1,x^2,x^3)\left[\begin{array}{ccc}-1& &\\
&1&\\&&-1\end{array}\right],$$
$$
\tau_{3}(x^1,x^2,x^3)=(x^1,x^2,x^3)\left[\begin{array}{ccc}-1& &\\
&-1&\\&&1\end{array}\right].$$ There is also an order $3$ automorphism $\sigma\in Aut(V_{L_{2}})$ defined
by
$$
\sigma(x^1,x^2,x^3)=(x^1,x^2,x^3)\left[\begin{array}{ccc}0&1 &0\\
0&0&-1\\-1&0&0\end{array}\right].$$ It is easy to see that $\sigma$ and
$\tau_{i},i=1,2,3$ generate a finite subgroup of $Aut(V_{L_{2}})$
isomorphic to the alternating group $A_{4}$. We simply denote this
subgroup by $A_{4}$. It is easy to check that the subgroup $K$
generated by $\tau_{i}$, $i=1,2,3$ is a normal subgroup of $A_{4}$
of order $4$. Let
$$J= h(-1)^4{\bf 1} -2h(-3)h(-1){\bf 1} + \frac{3}{2}h(-2)^2{\bf
1}, \ \ E=e^{\be}+e^{-\be}$$ where $h=\frac{1}{\sqrt{2}}\al$,
$\be=2\al$. The following lemma comes from \cite{DG} and \cite{DJ4}.
\begin{lem}\label{l3.1} (1) The vertex operator algebra $V_{L_{2}}^K$ and $V_{\Z\be}^{+}$ are the same and $V_{\Z\be}^{+}$ is generated by $J$ and $E$.
Moreover, $(V_{L_{2}}^K)_4$ is four dimensional with a basis consisting of
$L(-2)^2\1, L(-4)\1, J, E.$

(2)  The vertex operator algebra $V_{L_{2}}^{A_{4}}$ and $(V_{\Z\be}^+)^{\<\sigma\>}$ are the same.

(3) The action of $\sigma$ on $J$ and $E$ are given by
$$
\sigma(J)=-\frac{1}{2}J+\frac{9}{2}E, \ \
\sigma(E)=-\frac{1}{6}J-\frac{1}{2}E.
$$
\end{lem}

Clearly, $\sigma$ preserves the subspace of $(V_{\Z\beta}^+)_4$ spanned by $J$ and $E.$ It is easy to check that
\begin{equation}\label{e3.1}
\sigma(X^1)=\frac{-1+\sqrt{3}i}{2}X^1,
 \  \
\sigma(X^2)=\frac{-1-\sqrt{3}i}{2}X^2
\end{equation}
where
\begin{equation}\label{e3.9}X^1=J-\sqrt{27}iE, \ X^2=J+\sqrt{27}iE.
\end{equation}
This implies that $(V_{\Z\be}^+)^{\<\sigma\>}_4=L(1,0)_4$ where $L(1,0)$ is the vertex operator subalgebra
of $V_{\Z\beta}$ generated by $\omega.$ It follows from \cite{DG} that
$$(V_{\Z\be}^+)^{\<\sigma\>}=L(1,0)\bigoplus \sum_{n\geq 3}a_nL(1,n^2)$$
as a module for $L(1,0)$, where $a_n$ is the multiplicity of
$L(1,n^2)$ in $(V_{\Z\be}^+)^{\<\sigma\>}.$
Using (\ref{e3.1}) shows that for any $n\in\Z$,
$$
X^{1}_{n}X^2\in(V_{\Z\be}^+)^{\<\sigma\>}=V_{L_{2}}^{A_{4}}.
$$

We sometimes also call a highest weight vector for the Virasoro
algebra  a primary vector.  From \cite{DG} we know that $V_{\Z\be}^{+}$
contains two linearly independent primary vectors   $J$ and $E$ of
weight $4$ and one linearly independent primary vector of weight
$9.$ Note from \cite{DJ4} that
$$
J_{3}J=-72L(-4){\bf 1}+336L(-2)^2{\bf 1}-60J, \
E_{3}E=-\frac{8}{3}L(-4){\bf 1}+\frac{ 112}{9}L(-2)^2{\bf
1}+\frac{20}{9}J$$ (cf. \cite{DJ3}). By Theorem \ref{2t1} and Lemma
\ref{l3.1}, we have for $n\in\Z$
$$
X^1_{n}X^2\in L(1,0)\oplus L(1,9)\oplus L(1,16).
$$

The following lemma comes from \cite{DJ4}.
\begin{lem}\label{l3.4} The vector
\begin{eqnarray*}
& &u^{(9)}=-\frac{\sqrt{2}}{4}(J_{-2}E-E_{-2}J)\\
& &\ \ \ \ \  =-\dfrac{1}{\sqrt{2}}(15h(-4)h(-1)+10h(-3)h(-2)+10h(-2)h(-1)^3)\otimes
(e^{\be}+e^{-\be})\\
& &\ \ \ \ \ + (6h(-5)+10h(-3)h(-1)^2+\dfrac{15}{2}h(-2)^2h(-1)+h(-1)^5)\otimes
(e^{\be}-e^{-\be})
\end{eqnarray*}
is a non-zero primary vector of  $V_{\Z\be}^{+}$ of weight $9.$
\end{lem}

By Lemma \ref{l3.1}, we have $
J_{-9}J+27E_{-9}E\in (V_{\Z\be}^{+})^{\<\sigma\>}.$
Then
\begin{equation}\label{eq4.1}
J_{-9}J+27E_{-9}E=x^0+X^{(16)}+27(e^{2\be}+e^{-2\be}),
\end{equation}
where $x^0\in L(1,0)$, and  $X^{(16)}$  is a non-zero primary
element of weight 16 in $M(1)^{+}$. Denote
\begin{equation}\label{weight16}
u^{(16)}=X^{(16)}+27(e^{2\be}+e^{-2\be}).
\end{equation}
Then $u^{(16)}\in(V_{\Z\be}^{+})^{\<\sigma\>}$ is a non-zero primary
vector of weight 16. The following results come  from \cite{DJ4}.
\begin{theorem}\label{generators} The following hold:
(1) \ $V_{L_{2}}^{A_{4}}$ is generated by $u^{(9)}$.

(2) \  $V_{L_{2}}^{A_{4}}$ is linearly spanned by
$$
L(-m_{s})\cdots L(-m_{1})u^{(9)}_{n}u^{(9)}, \ L(-m_{s})\cdots
L(-m_{1})w^{p}_{-k_{p}}\cdots w^{1}_{-k_{1}}w,$$ where
$w,w^1,\cdots,w^p\in\{u^{(9)},u^{(16)}\}$, $k_{p}\geq \cdots\geq
k_{1}\geq 2$, $n\in\Z$, $m_{s}\geq\cdots \geq m_{1}\geq 1$, $s,p\geq
0$.
\end{theorem}
\begin{theorem}\label{rationality}
$V_{L_{2}}^{A_{4}}$ is $C_{2}$-cofinite and rational.
\end{theorem}

\section{Characterization of $V_{L_{2}}^{A_{4}}$}
\def\theequation{4.\arabic{equation}}
\setcounter{equation}{0}

In this section, we will give a characterization of the rational vertex operator algebra $V_{L_{2}}^{A_{4}}$. For this purpose we assume the following:

(A)  $V$ is a simple, $C_{2}$-cofinite£¬ rational£¬ CFT type  and self-dual vertex operator algebra of central charge 1;

(B) $V$ is a  sum of irreducible modules for the Virasoro algebra;

(C) There is a primary vector of weight $9$ and the weight of any  primary vector whose weight is greater than $0$ is greater than or equal to $9$.

%We first have the  following lemma from \cite{DJ1}.
%\begin{lem}\label{direct-sum}
% If $V$ is a simple vertex operator algebra such that $V_0=\C {\bf 1},$
 %$L(1)V_1=0$ and $c>1.$ Then $V$ is a completely reducible module for the
%Virasoro algebra.
%\end{lem}
%It follows from the  lemma \ref{direct-sum} that $V$ is a completely reducible module for the Virasoro algebra.

Obviously, $V_{L_{2}}^{A_{4}}$ satisfies (A)-(C). By Theorem \ref{dln} and Theorem \ref{Kir}, if a vertex operator algebra $V$ satisfies (A)-(C), then $V$ and $V_{L_{2}}^{A_{4}}$ have the same trace function.
So there is only one linearly independent primary vector of weight 9 in $V$.

For short, let  $V^{1}=V_{L_{2}}^{A_{4}}$ and $V^{2}$ be an arbitrary vertex operator algebra satisfying (A)-(C).  We will prove that $V^1$ and $V^2$ are isomorphic vertex operator algebras. Since  $V^1$ and $V^2$ have the same $q$-character, it follows from the assumption that $V^1$ and $V^2$ are isomorphic modules for the Virasoro algebra.
Let $x^i$ be a non-zero  weight 9 primary vectors in $V^i$, $i=1,2$ such that
\begin{equation}\label{e5.1}
(x^1,x^1)=(x^2,x^2).
\end{equation}
By \cite{DJ4}, there is only one linearly independent primary elements of weight 16 in $V^i$, $i=1,2$.  Now let $y^i\in V^i_{16}$ be  linearly independent primary vectors in $V^i$, $i=1,2$ such that
\begin{equation}\label{e5.2}
(y^1,y^1)=(y^2,y^2).        \ \              \ i=1,2.
\end{equation}

 The following lemma comes from  \cite{DJ4}.
 \begin{lem}\label{dj4.2}
  There is no non-zero primary vector of weight 25  in both $V^1$ and $V^2$.
  \end{lem}

Let $V^{(i,9)}$ be the $L(1,0)$-submodule of $V^i$ generated by $x^i$ and  $V^{(i,16)}$   the $L(1,0)$-submodule of $V^i$ generated by $y^i$, $i=1,2$. We may identify the Virasoro vertex operator subalgebra $L(1,0)$ both in $V^1$ and $V^2$. Let
 $$\phi:  \ L(1,0)\oplus V^{(1,9)}\oplus V^{(1,16)} \rightarrow L(1,0)\oplus V^{(2,9)}\oplus V^{(2,16)}$$
  be an $L(1,0)$-module isomorphism such that
  $$
  \phi\omega=\omega, \ \ \phi x^{1}=x^2, \ \ \phi y^1=y^2.$$
Then
$$
(u,v)=(\phi u,\phi v),
$$
for $u,v\in L(1,0)\oplus V^{(1,9)}\oplus V^{(1,16)}$.

Let
$$
{\cal I}^0(u,z)v={\cal P}^0\circ Y(u,z)v
$$
for $u,v\in V^{(1,9)}$ be the intertwining operator of type
$$\left(\hspace{-3 pt}\begin{array}{c}
L(1,0)\\ V^{(1,9)}\,V^{(1,9)}\end{array}\hspace{-3 pt}\right),$$
and
${\cal I}^0(\phi u,z)\phi v={\cal Q}^0\circ Y(\phi u,z)\phi v$ for $u,v\in V^{1}$ be the intertwining operator of type
$$\left(\hspace{-3 pt}\begin{array}{c}
L(1,0)\\ V^{(2,9)}\,V^{(2,9)}\end{array}\hspace{-3 pt}\right),$$
where ${\cal P}^{0}, {\cal Q}^0$ are the projections of $V^1$ and $V^2$ to $L(1,0)$ respectively. By (\ref{e5.1}), we have
\begin{equation}\label{e5.3}
{\cal I}^{0}(u,z)v={\cal I}^0(\phi u,z)\phi v,
\end{equation}
for $u,v\in V^{(1,9)}$.

Similarly, let $$
{\cal I}^1(u,z)v={\cal P}^1\circ Y(u,z)v
$$
for $u,v\in V^{(1,16)}$ be the intertwining operator of type
$$\left(\hspace{-3 pt}\begin{array}{c}
L(1,0)\\ V^{(1,16)}\,V^{(1,16)}\end{array}\hspace{-3 pt}\right),$$
and
${\cal I}^1(\phi u,z)\phi v={\cal Q}^1\circ Y(\phi u,z)\phi v$ for $u,v\in V^{1}$ be the intertwining operator of type
$$\left(\hspace{-3 pt}\begin{array}{c}
L(1,0)\\ V^{(2,16)}\,V^{(2,16)}\end{array}\hspace{-3 pt}\right),$$
where ${\cal P}^{1}, {\cal Q}^1$ are the projections of $V^1$ and $V^2$ to $L(1,0)$ respectively. By (\ref{e5.2}), we have
\begin{equation}\label{e5.4}
{\cal I}^{1}(u,z)v={\cal I}^1(\phi u,z)\phi v,
\end{equation}
for $u,v\in V^{(1,16)}$.

 Let
$$
{\cal I}^2(u,z)v={\cal P}^2\circ Y(u,z)v
$$
for $u\in V^{(1,16)}, v\in V^{(1,9)}$ be the intertwining operator of type
$$\left(\hspace{-3 pt}\begin{array}{c}
V^{(1,9)}\\ V^{(1,16)}\,V^{(1,9)}\end{array}\hspace{-3 pt}\right),$$
and
${\cal I}^2(\phi u,z)\phi v={\cal Q}^2\circ Y(\phi u,z)\phi v$ for $u\in V^{(1,16)}, v\in V^{(1,9)}$ be the intertwining operator of type
$$\left(\hspace{-3 pt}\begin{array}{c}
V^{(2,9)}\\ V^{(2,16)}\,V^{(2,9)}\end{array}\hspace{-3 pt}\right),$$
where ${\cal P}^{2}, {\cal Q}^2$ are the projections of $V^1$ and $V^2$ to $V^{(1,9)}$ and $V^{(2,9)}$ respectively.
Then  we have the following lemma.
\begin{lem}\label{l5.1}
Replacing $y^2$ by $-y^2$ if necessary, we have
$$
\phi({\cal I}^2(u,z)v)={\cal I}^2(\phi u,z)\phi v,
$$
for $u\in V^{(16)},v\in V^{(1,9)}$.
\end{lem}
\pf Since $V^{(1,9)}\cong V^{(2,9)} \cong L(1,9)$, $V^{(1,16)}\cong V^{(2,16)}\cong L(1,16)$, we may identify $V^{(1,9)}$ with $V^{(2,9)}$ and $V^{(1,16)}$ with $V^{(2,16)}$ through $\phi$. So both
${\cal I}^1(u,z)v$ and ${\cal I}^1(\phi u,z)\phi v$ for $u\in V^{(1,16)}, v\in V^{(9)}$ are intertwining operators of type
$$\left(\hspace{-3 pt}\begin{array}{c}
L(1,9)\\ L(1,16)\,L(1,9)\end{array}\hspace{-3 pt}\right).$$
Therefore
\begin{equation}\label{e5.5}
\phi({\cal I}^2(u,z)v)=\epsilon{\cal I}^2(\phi u,z)\phi v,
\end{equation}
for some $\epsilon\in\C$. By Theorem \ref{2t1} and (\ref{e5.1}), we have
$$
(y^1_{31}y^1)_{-1}x^1=(y^1,y^1)x^1, \ (y^2_{31}y^2)_{-1}x^2=(y^2,y^2)x^2.
$$
So
$$
\phi(y^1_{31}y^1)_{-1}x^1)=(y^2_{31}y^2)_{-1}x^2.
$$
On the other hand, we have
$$
(y^i_{31}y^i)_{-1}x^i=\sum\limits_{k=0}^{\infty}\left(\begin{array}{c}31\\k\end{array}\right)(-1)^k(y^i_{31-k}y^i_{-1+k}+y^i_{30-k}y^i_{k})x^i, \ i=1,2.
$$
Then by Theorem \ref{2t1} , Lemma \ref{dj4.2} and (\ref{e5.5}),
$$
(y^2_{31}y^2)_{-1}x^2=\epsilon^2((y^2_{31}y^2)_{-1}x^2).
$$
So we have $\epsilon^2=1$. If $\epsilon=1$, then the lemma holds. If $\epsilon=-1$, replacing $y^2$ by $-y^2$, then we get the lemma.\qed

Let
$$
{\cal I}^3(u,z)v={\cal P}^3\circ Y(u,z)v
$$
for $u,v\in V^{(1,9)}$ be the intertwining operator of type
$$\left(\hspace{-3 pt}\begin{array}{c}
V^{(1,16)}\\ V^{(1,9)}\,V^{(1,9)}\end{array}\hspace{-3 pt}\right),$$
and
${\cal I}^3(\phi u,z)\phi v={\cal Q}^3\circ Y(\phi u,z)\phi v$ for $u, v\in V^{(1,9)}$ be the intertwining operator of type
$$\left(\hspace{-3 pt}\begin{array}{c}
V^{(2,16)}\\ V^{(2,9)}\,V^{(2,9)}\end{array}\hspace{-3 pt}\right),$$
where ${\cal P}^{3}, {\cal Q}^3$ are the projections of $V^1$ and $V^2$ to $V^{(1,16)}$ and $V^{(2,16)}$ respectively.
Then  we have the following lemma.
\begin{lem}\label{l5.3}
$$
\phi({\cal I}^3(u,z)v)={\cal I}^3(\phi u,z)\phi v,
$$
for $u,v\in V^{(1,9)}$.
\end{lem}
\pf Note that both
${\cal I}^3(u,z)v$ and ${\cal I}^3(\phi u,z)\phi v$ for $u, v\in V^{(1,9)}$ are intertwining operators of type
$$\left(\hspace{-3 pt}\begin{array}{c}
L(1,16)\\ L(1,9)\,L(1,9)\end{array}\hspace{-3 pt}\right).$$
Therefore
\begin{equation}\label{e5.5'}
\phi({\cal I}^3(u,z)v)=\epsilon{\cal I}^3(\phi u,z)\phi v,
\end{equation}
for some $\epsilon\in\C$.  By Theorem \ref{2t1} and (\ref{e5.3}), we have
\begin{equation}\label{e5.6}
x^1_{1}x^1=u+a_{1}y^1, \ x^2_{1}x^2=u+a_2y^2
\end{equation}
where $u\in L(1,0)$ and $a_1,a_2\in \C$. By (\ref{e5.5'}),
\begin{equation}\label{e5.7}
a_{1}=\epsilon a_{2}.
\end{equation} Then by Theorem \ref{2t1}, we have
\begin{equation}\label{e5.8}x^1_{0}x^1=v+a_{1}L(-1)y^1, \ x^2_{0}x^2=v+a_2L(-1)y^2,
\end{equation}
for some $v\in L(1,0)$. Notice that
$$
(x^1_{1}x^1,y^1)=a_{1}(y^1,y^1), \ (x^2_{1}x^2,y^2)=a_{2}(y^2,y^2).$$
By (\ref{e5.5})-(\ref{e5.8}),
\begin{equation}\label{e5.9}(x^1_{1}x^1,y^1)=\epsilon (x^2_{1}x^2,y^2).
\end{equation}
On the other hand, we have
$$
(x^1_{1}x^1,y^1)=-(x^1,x^1_{15}y^1)=-(x^1,y^1_{15}x^1), \ (x^2_{1}x^2,y^2)=-(x^2,x^2_{15}y^2)=-(x^2,y^2_{15}x^2).
$$
By Lemma \ref{l5.3},
$$
\phi(y^1_{15}x^1)=y^2_{15}x^2.$$ So
$$(x^1_{1}x^1,y^1)=(x^2_{1}x^2,y^2).
$$
This together with (\ref{e5.9}) deduces that $\epsilon=1$. \qed

Let
$$
{\cal I}^4(u,z)v={\cal P}^4\circ Y(u,z)v
$$
for $u,v\in V^{(1,16)}$ be the intertwining operator of type
$$\left(\hspace{-3 pt}\begin{array}{c}
V^{(1,16)}\\ V^{(1,16)}\,V^{(1,16)}\end{array}\hspace{-3 pt}\right),$$
and
${\cal I}^4(\phi u,z)\phi v={\cal Q}^4\circ Y(\phi u,z)\phi v$ for $u, v\in V^{(1,16)}$ be the intertwining operator of type
$$\left(\hspace{-3 pt}\begin{array}{c}
V^{(2,16)}\\ V^{(2,16)}\,V^{(2,16)}\end{array}\hspace{-3 pt}\right),$$
where ${\cal P}^{4}, {\cal Q}^4$ are the projections of $V^1$ and $V^2$ to $V^{(1,16)}$ and $V^{(2,16)}$ respectively.
Then  we have the following lemma.
\begin{lem}\label{l5.6}
$$
\phi({\cal I}^4(u,z)v)={\cal I}^4(\phi u,z)\phi v,
$$
for $u,v\in V^{(1,16)}$.
\end{lem}
\pf It suffices to prove that
$$
\phi(y^1_{15}y^1)=y^2_{15}y^2.$$
By Lemma \ref{l5.3} and (\ref{e5.6}), we have
\begin{equation}\label{e5.10}
\phi(x^1_{1}x^1)=\phi(u+a_{1}y^1)=u+a_{1}y^2,
\end{equation}
where $u$ and $a_{1}$ are as in (\ref{e5.6}). Notice that
$$
(x^i_{1}x^i)_{15}y^i=\sum\limits_{k=0}^{\infty}\left(\begin{array}{c}1\\k\end{array}\right)(-1)^k(x^i_{1-k}x^i_{15+k}+x^i_{16-k}x^i_{k})y^i, \ i=1,2.
$$
Then by Lemma \ref{l5.3}, Lemma \ref{l5.1}, and the skew-symmetry property of vertex operator algebras, we have
$$
\phi((x^1_{1}x^1)_{15}y^1)=(x^2_{1}x^2)_{15}y^2.
$$
This together with (\ref{e5.10}) deduces that $\phi(y^1_{15}y^1)=y^2_{15}y^2$. The lemma follows.
\qed

Summarizing  lemmas \ref{l5.1}-\ref{l5.6}, we have the following proposition:
\begin{prop}\label{lfu5.4} (1) For any $u^1,v^1\in L(1,0)\oplus V^{(1,9)}\oplus V^{(1,16)}$, we have
$$
\phi(u^1_{n}v^1)=(\phi u^1)_{n}(\phi v^1),
$$
for $n\in\N$.

(2) For any $u^1,v^1\in L(1,0)\oplus V^{(1,9)}\oplus V^{(1,16)}$, we have
$$
(u^1,v^1)=(\phi(u^1),\phi(v^1)).
$$
\end{prop}

Recall from  Theorem \ref{generators} that $V^1$ is generated by $x^1$ and $V^1$ is linearly spanned by
$$
L(-m_{s})\cdots L(-m_{1})x^1_{n}x^1, \ L(-m_{s})\cdots
L(-m_{1})u^{p}_{-k_{p}}\cdots u^{1}_{-k_{1}}v,$$ where $x^1,y^1$ are the same as above and
$v,u^1,\cdots,u^p\in\{x^1,y^1\}$, $k_{p}\geq \cdots\geq
k_{1}\geq 2$, $n\in\Z$, $m_{s}\geq\cdots \geq m_{1}\geq 1$, $s,p\geq
0$. Our goal next is to show that $V^2$ is generated by  $\phi(x^1)=x^2$ and
 has a similar spanning set.

\begin{lem}\label{l5.8} For any $k,l\geq 1, s_{i},t_{i},p_{i}\geq 0$, $m_{i1},\cdots,m_{is_{i}}, n_{j1},\cdots,n_{jt_{j}}, r_{j1},\cdots, r_{jp_{j}}\in\Z_{+}$, $n_j\in\Z$, $u^{j1},\cdots,u^{jp_{j}},u^j\in\{x^1,y^1\}$, $i=1,2,\cdots,k,j=1,2,\cdots,l$,   if
\begin{eqnarray*}
& u'=\sum\limits_{i=1}^{k}a_{i}L(-m_{i1})\cdots L(-m_{is_{i}})x^2_{n_{i}}x^2\\
& +\sum\limits_{j=1}^lb_{j}L(-n_{j1})\cdots L(-n_{jt_{j}})(\phi u^{j1})_{-r_{j1}}
\cdots (\phi u^{jp_{j}})_{-r_{jp_{j}}}(\phi u^j)=0
\end{eqnarray*}
for some $a_i,b_j\in\C$
 then
\begin{eqnarray*}
& u=\sum\limits_{i=1}^{k}a_{i}L(-m_{i1})\cdots L(-m_{is_{i}})x^1_{n_{i}}x^1\\
& +\sum\limits_{j=1}^lb_{j}L(-n_{j1})\cdots L(-n_{jt_{j}})u^{j1}_{-r_{j1}}
u^{j2}_{-r_{j2}}\cdots u^{jp_{j}}_{-r_{jp_{j}}}u^j=0.
\end{eqnarray*}
\end{lem}
\pf  Without loss, we may assume that $u$ is a linear combination of homogeneous elements with same weight. Suppose that $u\ne 0.$ Since $V^{1}$ is self-dual and generated by $x^1,$
there is $x^1_{r_{1}}x^1_{r_{2}}\cdots x^1_{r_{q}}x^1\in V^1$ such that
\begin{equation}\label{bili-3}
(u, x^1_{r_{1}}x^1_{r_{2}}\cdots x^1_{r_{q}}x^1)\neq 0.
\end{equation}
{\bf Claim:} \ For any $v,u^1,\cdots,u^p\in\{x^1,y^1\}$, $k_{p}\geq \cdots\geq
k_{1}\geq 2$, $q_1,q_2,\cdots,q_t, n\in\Z$, $m_{s}\geq\cdots \geq m_{1}\geq 1$, $s,p,t\geq
0$,
\begin{eqnarray}
& &(L(-m_{s})\cdots L(-m_{1})x^1_{n}x^1, x^1_{q_1}x^1_{q_2}\cdots x^1_{q_t}x^1)\nonumber\\
& & \ \ \ \ =(L(-m_{s})\cdots L(-m_{1})x^2_{n}x^2, x^2_{q_1}x^2_{q_2}\cdots x^2_{q_t}x^2),\label{bili-1}
\end{eqnarray}
\begin{eqnarray}
& &(L(-m_{s})\cdots L(-m_{1})u^{p}_{-k_{p}}\cdots u^{1}_{-k_{1}}v, x^1_{q_1}x^1_{q_2}\cdots x^1_{q_t}x^1)\nonumber\\
& &\ \ \ \ =(L(-m_{s})\cdots L(-m_{1})\phi(u^{p})_{-k_{p}}\cdots \phi(u^{1})_{-k_{1}}\phi(v), x^2_{q_1}x^2_{q_2}\cdots x^2_{q_t}x^2).\label{bili-2}
\end{eqnarray}

We only show (\ref{bili-2}) as the proof for (\ref{bili-1}) is similar and simpler. We may assume that
$$\wt(L(-m_{s})\cdots L(-m_{1})u^{p}_{-k_{p}}\cdots u^{1}_{-k_{1}}v)=\wt( x^1_{q_1}x^1_{q_2}\cdots x^1_{q_t}x^1).$$ We prove (\ref{bili-2}) by induction on $\wt(L(-m_{s})\cdots L(-m_{1})u^{p}_{-k_{p}}\cdots u^{1}_{-k_{1}}v)$. By Proposition \ref{lfu5.4}, (\ref{bili-1}) holds if $\wt(L(-m_{s})\cdots L(-m_{1})u^{p}_{-k_{p}}\cdots u^{1}_{-k_{1}}v)<36$. If $s\geq 1$, then
\begin{eqnarray*}
& &(L(-m_{s})\cdots L(-m_{1})u^{p}_{-k_{p}}\cdots u^{1}_{-k_{1}}v, x^1_{q_1}x^1_{q_2}\cdots x^1_{q_t}x^1)\\
& &\ \ \ \ =(L(-m_{s-1})\cdots L(-m_{1})u^{p}_{-k_{p}}\cdots u^{1}_{-k_{1}}v, L(m_{s})x^1_{q_1}x^1_{q_2}\cdots x^1_{q_t}x^1),\\
& &(L(-m_{s})\cdots L(-m_{1})\phi(u^{p})_{-k_{p}}\cdots \phi(u^{1})_{-k_{1}}\phi(v), x^2_{q_1}x^2_{q_2}\cdots x^2_{q_t}x^2)\\
& &\ \ \ \ =(L(-m_{s-1})\cdots L(-m_{1})\phi(u^{p})_{-k_{p}}\cdots \phi(u^{1})_{-k_{1}}\phi(v), L(m_{s})x^2_{q_1}x^2_{q_2}\cdots x^2_{q_t}x^2).
\end{eqnarray*}
So by inductive assumption, (\ref{bili-2}) holds.

If $s=0$, then
\begin{eqnarray*}
& &(u^{p}_{-k_{p}}\cdots u^{1}_{-k_{1}}v, x^1_{q_1}x^1_{q_2}\cdots x^1_{q_t}x^1)\\
& &\ \ \ \ =(u^{p-1}_{-k_{p-1}}\cdots u^{1}_{-k_{1}}v, u^p_{2\wt(u^p)+k_{p}-2}x^1_{q_1}x^1_{q_2}\cdots x^1_{q_t}x^1),
\end{eqnarray*}
\begin{eqnarray*}
& &(\phi(u^{p})_{-k_{p}}\cdots \phi(u^{1})_{-k_{1}}\phi(v), x^2_{q_1}\cdots x^2_{q_t}x^2)\\
& &\ \ \ \ =(\phi(u^{p-1})_{-k_{p-1}}\cdots \phi(u^{1})_{-k_{1}}\phi(v), \phi(u^p)_{2\wt(u^p)+k_{p}-2}x^2_{q_1}\cdots x^2_{q_t}x^2).
\end{eqnarray*}
Since $k_{p}\geq 2$, by inductive assumption, (\ref{bili-2}) holds if $u^p=x^1.$  If $u^p=y^1$ by (\ref{e5.6})
$u^p_{2\wt(u^p)+k_{p}-2}$ is a sum of operators of forms
$aL(n_1)\cdots L(n_\mu), bx^1_ix^1_j$  of the same weight where $n_1\leq \cdots \leq n_{\mu}$ and all $n_t$ are nonzero.  By induction assumption we know that
\begin{eqnarray*}
& &(u^{p-1}_{-k_{p-1}}\cdots u^{1}_{-k_{1}}v, x^1_ix^1_jx^1_{q_1}x^1_{q_2}\cdots x^1_{q_t}x^1)\\
& &\ \ \ \ =(\phi(u^{p-1})_{-k_{p-1}}\cdots \phi(u^{1})_{-k_{1}}\phi(v), x^2_ix^2_jx^2_{q_1}\cdots x^2_{q_t}x^2).
\end{eqnarray*}
Also by  Proposition \ref{lfu5.4}, relation (\ref{e5.10}), the fact that $x^i$ are highest weight vectors for the Virasoro algebra with the same weight, and the invariant properties of the bilinear forms,
\begin{eqnarray*}
& &(u^{p-1}_{-k_{p-1}}\cdots u^{1}_{-k_{1}}v, L(n_1)\cdots L(n_\mu)x^1_{q_1}x^1_{q_2}\cdots x^1_{q_t}x^1)\\
& &\ \ \ \ =(\phi(u^{p-1})_{-k_{p-1}}\cdots \phi(u^{1})_{-k_{1}}\phi(v), L(n_1)\cdots L(n_\mu)x^2_{q_1}\cdots x^2_{q_t}x^2).
\end{eqnarray*}
So the claim is proved.

By the claim and (\ref{bili-3}), we have
$$
(u', x^1_{r_{1}}x^2_{r_{2}}\cdots x^2_{r_{q}}x^2)\neq 0,
$$
which contradicts the assumption that $u'=0$. \qed

Let $U^2$ be the subalgebra of $V^2$ generated by $x^2$ and $y^2$. By Theorem \ref{generators} and Lemma \ref{l5.8}, for every $n\geq 0$, $\dim V^1_{n}\leq \dim U^2_n$. Since $V^1$ and $V^2$ have the same
 graded dimensions, it follows that $\dim V^1_n=\dim V^2_n$ for $n\geq 0$. So $\dim V^2_n=\dim U^2_n$ for $n\geq 0$ and $V^2=U^2$. So we have the following corollary which is essentially  the $V^2$ version of Theorem \ref{generators}.

\begin{coro}\label{corol-1}
 $V^2$ is linearly spanned by
$$
L(-m_{s})\cdots L(-m_{1})x^2_{n}x^2, \ L(-m_{s})\cdots
L(-m_{1})v^{p}_{-k_{p}}\cdots v^{1}_{-k_{1}}v,$$ where $x^2,y^2$ are the same as above and
$v,v^1,\cdots,v^p\in\{x^2,y^2\}$, $k_{p}\geq \cdots\geq
k_{1}\geq 2$, $n\in\Z$, $m_{s}\geq\cdots \geq m_{1}\geq 1$, $s,p\geq
0$.
\end{coro}
Define $\psi(x^2)=x^1, \psi(y^2)=y^1$, and extend $\psi$ to $\psi: V^2\rightarrow V^1$ by
$$
\psi(L(-m_{s})\cdots L(-m_{1})x^2_{n}x^2)=L(-m_{s})\cdots L(-m_{1})x^1_{n}x^1
$$
and
$$
L(-m_{1})v^{p}_{-k_{p}}\cdots v^{1}_{-k_{1}}v=L(-m_{1})\psi(v^{p})_{-k_{p}}\cdots \psi(v^{1})_{-k_{1}}\psi(v),
$$
where $v,v^1,\cdots,v^p\in\{x^2,y^2\}$, $k_{p}\geq \cdots\geq
k_{1}\geq 2$, $n\in\Z$, $m_{s}\geq\cdots \geq m_{1}\geq 1$, $s,p\geq
0$. Then by the  discussion above, $\psi$ is a linear isomorphism from $V^2$ to $V^1$.  It follows that $\phi$ can be extended to a linear isomorphism from $V^1$ to $V^2$ such that
$$
\phi(L(-m_{s})\cdots L(-m_{1})x^1_{n}x^1)=L(-m_{s})\cdots L(-m_{1})x^2_{n}x^2
$$
and
$$
L(-m_{1})u^{p}_{-k_{p}}\cdots u^{1}_{-k_{1}}u=L(-m_{1})\phi(p^{p})_{-k_{p}}\cdots \phi(u^{1})_{-k_{1}}\phi(u),
$$
where $u,u^1,\cdots,u^p\in\{x^1,y^1\}$, $k_{p}\geq \cdots\geq
k_{1}\geq 2$, $n\in\Z$, $m_{s}\geq\cdots \geq m_{1}\geq 1$, $s,p\geq
0$.

We are now in a position to state our main result of this paper.
\begin{theorem}\label{main-theorem}
If a vertex operator algebra $V$ satisfies the conditions (A)-(C), then   $V$ is isomorphic to $V_{L_{2}}^{A_{4}}$.
\end{theorem}
\pf Recall that $V^1\cong V_{L_{2}}^{A_{4}}$ satisfying (A)-(C). So it suffices to show that $\phi$ is a vertex operator algebra automorphism from $V^1$ to $V^2$. Let $u=x^1_{m_1}x^1_{m_2}\cdots x^1_{m_s}x^1$, $v=x^1_{q_1}x^1_{q_2}\cdots x^1_{q_t}x^1\in V^1$, where $m_i,q_j\in \Z, i=1,2,\cdots,p_s,j=1,2,\cdots,t$. We need to show that for any $n\in\Z$, $\phi(u^1_nu^2)=\phi(u^1)_n\phi(u^2)$. Note from Theorem \ref{2t1} that for $m_{1},m_{2}\in\Z_{+}$,  $x^1_{m_{1}}x^1\in L(1,0)\oplus V^{(1,16)}$, $y^1_{m_{2}}x^1\in V^{(1,9)}$. Since  for any $p,q\in\Z$,
$$
x^i_{q}x^i_{p}=x^i_px^i_q+\sum\limits_{j=0}^{\infty}\left(\begin{array}{c}q\\j\end{array}\right)(x^i_jx^i)p+q-j, \ i=1,2,
$$
$$
y^i_{q}x^i_{p}=y^i_px^i_q+\sum\limits_{j=0}^{\infty}\left(\begin{array}{c}q\\j\end{array}\right)(y^i_jx^i)p+q-j,
$$
Then by Lemma \ref{lfu5.4}, it is easy to see that for any fixed $n\in\Z$,
\begin{eqnarray*}
& u^1_nu^2=\sum\limits_{i=1}^{k}a_{i}L(-m_{i1})\cdots L(-m_{is_{i}})x^1_{n_{i}}x^1\\
& +\sum\limits_{i=1}^lb_{i}L(-n_{i1})\cdots L(-n_{it_{i}})u^{i1}_{-r_{i1}}
u^{i2}_{-r_{i2}}\cdots u^{ip_{i}}_{-r_{ip_{i}}}u^i,
\end{eqnarray*}
for some $k,l\geq 1, s_{i},t_{i},p_{i}\geq 0$, $m_{i1},\cdots,m_{is_{i}}, n_{j1},\cdots,n_{jt_{j}}, r_{j1},\cdots, r_{jp_{j}}\in\Z_{+}$, $n_i\in\Z$, $u^{i1},\cdots,u^{ip_{i}},u^i\in\{x^1,y^1\}$, $i=1,2,\cdots,k,j=1,2,\cdots,l$, then
\begin{eqnarray*}
& \phi(u^1)_n\phi(u^2)=\sum\limits_{i=1}^{k}a_{i}L(-m_{i1})\cdots L(-m_{is_{i}})x^2_{n_{i}}x^2\\
& +\sum\limits_{i=1}^lb_{i}L(-n_{i1})\cdots L(-n_{it_{i}})(\phi u^{i1})_{-r_{i1}}
\cdots (\phi u^{ip_{i}})_{-r_{ip_{i}}}(\phi u^i).
\end{eqnarray*}
The proof is complete. \qed

{\bf Acknowledgments}  The first author acknowledges the support  from NSF grants and a
Faculty research grant from  the University of California at Santa
Cruz. The second author acknowledges the support from China NSF grants
(10931006 and 11371245), the RFDP grants of China (20100073110052), and the Innovation
Program of Shanghai Municipal
Education Commission (11ZZ18).

\end{document}